\documentclass[12pt]{article}
 \usepackage{amsmath}
\usepackage{amsfonts}
\usepackage{tikz}
\usepackage[title]{appendix}
\usepackage{mathrsfs}
\usepackage{subfig}
\usepackage{array}
\usepackage{blkarray}
\usepackage{tabu}
\usepackage{nth}
\usepackage[
  separate-uncertainty = true,
  multi-part-units = repeat
]{siunitx}

\usetikzlibrary{arrows}
\usepackage{epsfig,}
\usepackage{epsfig}
\tikzset{
    vertex/.style = {
        circle,
        draw,
        outer sep = 3pt,
        inner sep = 3pt,
    },edge/.style = {->,> = latex'}
}

\usepackage{amssymb}
\usepackage{enumitem}
\usepackage{amsthm}
\usepackage{dsfont}
\def\diag{\mathop{\rm diag}}

\def\Diag{\mathop{\rm Diag}}

\def\rank{\mathop{\rm rank}}
\newcommand{\rr}{\mathbb{R}}
\newcommand{\0}{\mathbf{0}}

\newtheorem{theorem}{Theorem}

\newtheorem{corollary}{Corollary}

\newtheorem{example}{Example}

\newtheorem{lemma}{Lemma}
\newtheorem{definition}{Definition}

\begin{document}
\begin{center}
\begin{large}
Conditionally Euclidean distance matrices
\end{large}
\end{center}
\begin{center}
Shivani Goel and Hiroshi Kurata\\
\today
\end{center}

\begin{abstract}
We introduce a new class of matrices called \textit{conditional Euclidean distance matrices} in this paper. These matrices are a generalisation of Euclidean distance matrices. We develop the theory of conditional Euclidean distance matrices and show that these matrices represent the squared Lorentzian distance between points of a  Minkowski space. We also highlight the connection of  conditional Euclidean distance matrices with representation of hierarchical data in machine learning and the theory of general relativity in physics. 
\end{abstract}
{\bf AMS CLASSIFICATION.} 15A48, 15B48, 05C50

\section{Introduction}
Let $\mathscr{S}^n$  denote the set of all $n \times n$ symmetric matrices. A hollow symmetric matrix is a symmetric matrix with zero diagonal elements. Let $\mathscr{S}_h^n$ be
the set of hollow symmetric matrices of order $n$. By a centered symmetric matrix we mean a symmetric matrix with zero row (and hence column) sums. Let $\mathscr{S}_{e}^n$ be the set of all centered symmetric matrices of order $n$ given by:
\[\mathscr{S}_{e}^n:=\{F \in \mathscr{S}^n: Fe=\0\}.\]
 As shown by Critchley \cite{CRITCHLEY198891}, Theorem 2.2, and
Johnson and Tarazaga \cite{JOHNSON1995375}, pp.380, there is a one-to-one correspondence between $\mathscr{S}_h^n$ and $\mathscr{S}_{e}^n$. The following two linear maps $\mathscr{T}$ and $\mathscr{K}$ define this one-to-one correspondence:
\[\mathscr{T}: \mathscr{S}_h^n \rightarrow \mathscr{S}_{e}^n: D \rightarrow -\frac{1}{2}(I-\frac{ee^T}{n})D(I-\frac{ee^T}{n}),\]
\[\mathscr{K}: \mathscr{S}_{e}^n \rightarrow \mathscr{S}_h^n: F \rightarrow \Diag(F)E+E\Diag(F)^T-2F.\]
Let $\mathscr{D}^n$ denote the set of $n \times n$ Euclidean distance matrices (see (\ref{Def_EDM})) and let 
\[\Omega_n(e):= \{F \in \mathscr{S}_{e}^n: F~\mbox{is positive semidefinite}\}. \]
It is of interest to note that both $\mathscr{D}^n$ and $\Omega_n(e)$ forms convex cones in $\mathscr{S}_h^n$ and $\mathscr{S}_{e}^n$, respectively. It is well known that the mappings $\mathscr{T}|_{\mathscr{D}^n}$ and $\mathscr{K}|_{\Omega_n(e)}$ are mutually inverse. Consider the following subset of $\mathscr{S}_{e}^n$:
\[\Omega_n^-(e):= \{F \in \mathscr{S}_{e}^n: F~\mbox{has at most one negative eigenvalue}\}. \]
Clearly $\Omega_n(e) \subset \Omega_n^-(e)$. Further, $\Omega_n^-(e)$ forms a cone in $\mathscr{S}_{e}^n$. In this paper we are interested in identifying and studying the image of $\Omega_n^-(e)$ under the linear map $\mathscr{K}$. 

\subsection{Historical background}
The motivation for work on hollow symmetric matrices with non-negative off diagonal entries takes its root from Euclidean geometry and originated in the classical paper \cite{schoenberg1935remarks} of Schoenberg. Euclidean distance matrices were first studied formally in the work \cite{gower1966}, \cite{gower1982} and \cite{GOWER1985} of Gower. Euclidean distance matrices were characterized into their modern definition by Hayden and Wells in \cite{HAYDEN1988115}. 

An $n \times n$  hollow symmetric matrix $D$ is called an \textit{Euclidean distance matrix} (EDM) if for some positive integer $r$, there exist points $p^1, \dotsc , p^n \in \rr^r$ such that
\begin{equation}\label{Def_EDM}
    d_{ij}= \|p^i-p^j\|^2
\end{equation}
for each $i, j \in \{1,\dotsc,n\}$. Several characterizations and results on EDMs can be found in the works of Tarazaga et al. \cite{TARAZAGA199677} and \cite{Tarazaga2007}, Hayden et al.  \cite{HAYDEN199997}, Tarazaga \cite{TARAZAGA20051}, Alfakih \cite{ALFAKIH2006348}, Bénasséni \cite{BENASSENI2006365}, Kurata et al. \cite{KURATA2007586}, \cite{KURATA2010534}, \cite{KURATA20133177} and \cite{KURATA2015106}, Balaji et al. \cite{BALAJIedm} and many more references within these articles. Kurata and Bapat in \cite{KurataBapat+2016} considered the class of all hollow symmetric matrices and extended some of the properties of EDMs to this class. Their work motivates us to look at a subset of the class of hollow symmetric matrices that properly contains the class of EDMs, but most of the results on EDMs can be naturally generalised to this new class.  

\subsection{Notations and preliminaries}
In this section we introduce some notations and preliminaries relevant to this work. We omit including the notations and preliminaries defined in the previous sections to avoid redundancy. Further some minimally used notations will be introduced as and when needed to 
abstain inundating the reader.  

The notations $e$ and $\0$ are used to represent the vectors of relevant order with all entries $1$ and $0$, respectively. The identity and all ones matrices are represented by $I$ and $E$, respectively. The transpose of a vector $x$ is denoted by $x^T$, and we use $A^T$ to represent the transpose of a matrix. Two vectors $x$ and $y$ are orthogonal if $x^Ty=0$. For any vector $x$ and subspace $W$, $x^\perp$ and $W^\perp$ denote the orthogonal complement of ${\rm span}\{x\}$ and $W$, respectively. The notation $A=(a_{ij})$ is used to represent a matrix with its $ij^{\rm th}$ entry as $a_{ij}$. The rank of a matrix $A$ is denoted by $\rank(A)$. For a matrix $A$, $\Diag(A)$ represent the diagonal matrix formed by the diagonal entries of $A$. The matrix $A^\dag$ is the Moore-Penrose of $A$. 

A matrix $A \in \mathscr{S}^n$ is called \textit{positive semidefinite} if and only if $x^TAx \geq 0$ for all $x \in \rr^n$. A famous characterization for positive semidefinite matrices is that a symmetric matrix is positive semidefinite if and only if all of its eigenvalues are non-negative. A matrix $A$ is said to be positive (negative) semidefinite on a subspace $W$, if $x^TAx \geq 0$ ($x^TAx \leq 0$) for any $x \in W$. 

\subsection{Outline of the paper}
In section 2, we introduce the class of conditional Euclidean matrices and explore its characterisation from various perspectives. We also study properties of the centered symmetric matrix associated with a coEDM in section 2. Next, the concept of Gale space is generalised for a  hollow symmetric matrix in section 3 and we obtain similar results about these spaces as are known for EDMs. In section 4, we discuss some special properties of non-negative coEDMs. We conclude this paper by showing a connection of coEDMs with Minowski space, machine learning and in the theory of general relativity in section 5.

\section{Conditional Euclidean distance matrix}
For the rest of the paper, we use $J$ to denote the matrix $I-\frac{1}{n}E$. Suppose $D$ is a 
hollow-symmetric matrix and let $F=-\frac{1}{2}JDJ$. It is well-known from \cite{schoenberg1935remarks} that $D$ is an EDM if and only if $F$ is a positive semidefinite matrix. We are interested in 
hollow-symmetric matrices for which $F$ is not positive semidefinite. A matrix that is not positive semidefinite has at least one negative eigenvalue. We recall  
\[\Omega_n^-(e)= \{F \in \mathscr{S}_{e}^n: F~\mbox{has at most one negative eigenvalue}\}, \]
and define the following.

\begin{definition}\label{new_def_coEDM}\rm
    A 
     hollow symmetric matrix $D$ is a \textit{conditional Euclidean distance matrix} (coEDM) if and only $D \in \mathscr{K}(\Omega_n^-(e))$.  We denote the set of all $n \times n$ coEDMs by $\mathscr{C}^n$.
\end{definition}

Since the image of a cone under a linear map is a cone, we conclude that the set of all $n \times n$ coEDMs forms a cone in $\mathscr{S}_h^n$. Let 
\[\Omega_n^{(1)}(e) := \{F \in \mathscr{S}_e^n | F \mbox{ has exactly one negative eigenvalue} \}.\]
Then it is clear that the set $\Omega_n^-(e)$ can be decomposed as the following disjoint union 
\[\Omega_n^-(e)= \Omega_n(e) \cup \Omega_n^{(1)}(e).\] 
Since $\mathscr{K}(\Omega_n^-(e)) = \mathscr{C}^n$, we now have
\[ \mathscr{C}^n= \mathscr{K}(\Omega_n^-(e)) = \mathscr{K}(\Omega_n(e)) \cup \mathscr{K}(\Omega_n^{(1)}(e)). \]
Since the set $\mathscr{D}^n$ of $n \times n$ EDMs is given by $ \mathscr{K}(\Omega_n(e))$, we have
\begin{equation}\label{coEDM_is_union}
 \mathscr{C}^n=\mathscr{D}^n \cup \mathscr{K}(\Omega_n^{(1)}(e)). 
\end{equation}
Thus,  every EDM is a coEDM. However, the converse is not true. The following example shows such an instance.
\begin{example} \rm
    Consider the matrix 
    \[D= \frac{1}{4}\left[\begin{array}{ccc}
0 & -19 & 8 \\
-19 & 0 & 41 \\
8 & 41 & 0 
\end{array}\right].\]
Clearly, $D$ is not an EDM due to having negative entries. But it is a coEDM since the matrix $D$ can be written as $D=\mathscr{K}(F)$, where
\[F=-\frac{1}{2}JDJ = \frac{1}{4}\left[\begin{array}{rrr}
-7 & 8 & -1 \\
8 & 4 & -12 \\
-1 & -12 & 13
\end{array}\right]\]
is a centered symmetric matrix having exactly one negative eigenvalue. 
\end{example}

 In the next section, we study $\mathscr{C}^n$, the set of all coEDMs in detail and provide some characterizations. 

\subsection{Characterization of coEDMs}
Suppose $D$ is a 
hollow-symmetric matrix. From the results of Critchley \cite{CRITCHLEY198891} and Johnson \cite{JOHNSON1995375}, there exists a centered symmetric matrix $F$ such that 
\begin{equation*} \label{hollowsymmetric}
    D=\mathscr{K}(F)=E\Diag(F)+\Diag(F)E-2F.
\end{equation*} 
Since $F$ has zero row sums, it is easy to see that
\begin{equation}\label{form2}
    F=-\frac{1}{2}(I-es^T)D(I-se^T),
 \end{equation}
where $s=\frac{1}{n}e$. This shows that if $ D=\mathscr{K}(F)$, then $F$ is uniquely determined. Thus, if $D$ is a coEDM, then by Theorem \ref{new_def_coEDM}, $F \in \Omega_n^-(e)$. Since $F$ has at most one negative eigenvalue, there exists $\omega$ such that $F$ is positive semidefinite on $\omega^\perp$. Since $F$ is a zero row sum matrix, we can assume that $\omega^Te=0$.  Using the above discussion, we obtain the following characterization of coEDMs.

\begin{theorem}\label{def_coEDM}
     A 
     hollow symmetric matrix $D$ is a \textit{conditional Euclidean distance matrix} (coEDM) if and only if the matrix $F=-\frac{1}{2}JDJ$ is positive semidefinite on $\omega^\perp$ for some $\omega \in e^\perp$. 
\end{theorem}

Gower in \cite{gower1982} showed that $D$ is a Euclidean distance matrix if and only if (\ref{form2}) is positive semidefinite for any $s$ with $s^Te=1$ and $s^TD \neq 0$. Later, Gower in \cite{GOWER1985} showed that the condition $s^TD \neq 0$ is superfluous. Thus $D$ is a Euclidean distance matrix if and only if (\ref{form2}) is positive semidefinite for any $s$ with $s^Te=1$. The next result provides a generalisation of Gower's result for coEDMs. 

\begin{theorem}\label{proj_form}
    Let $D$ be a  
     hollow symmetric $n \times n$ matrix. Then $D$ is a conditionally Euclidean distance matrix if and only if there exist $\omega \in e^\perp$ such that the matrix 
    \[ -\frac{1}{2}(I-es^T)D(I-se^T)\]
    is positive semidefinite on $\omega^\perp$ for any $s$ such that $s^Te=1$, and $s^T\omega=0$.
\end{theorem}
\begin{proof}
    Suppose $\omega \in e^\perp$ such that
    \[L = -\frac{1}{2}(I-es^T)D(I-se^T)\]
    is positive semidefinite on $\omega^\perp$ for any $s$ such that $s^Te=1$, and $s^T\omega=0$. Since $\frac{1}{n}e^Te = 1$, $\frac{1}{n}e^T\omega = 0$, $-\frac{1}{2}JDJ$ is positive semidefinite on $\omega^\perp$
and hence $D$ is a coEDM. 
    
    Conversely, suppose $D$ is a coEDM and let $s=\frac{1}{n}e$. Then there exists $\omega \in e^\perp$ such that
    \[F = -\frac{1}{2}(I-es^T)D(I-se^T)\]
    is positive semidefinite on $\omega^\perp$. Consider a vector $t$ such that $t^Te=1$ and $t^T\omega=0$. Then for $s=\frac{1}{n}e$,
\[(I-et^T)(I-es^T) = I-et^T,\]
and 
\[(I-et^T)F(I-te^T) = -\frac{1}{2}(I-et^T)D(I-te^T).\]
Suppose $x \in \omega^\perp$, then 
\begin{equation*}
    \begin{aligned}
         -\frac{1}{2}x^T(I-et^T)D(I-te^T)x &=x^T(I-et^T)F(I-te^T)x \geq 0.
    \end{aligned}
\end{equation*}
Thus, $-\frac{1}{2}(I-es^T)D(I-se^T)$ is positive semidefinite on $\omega^\perp$ for any $s$ such that $s^Te=1$ and $s^T\omega=0$. The proof is complete.
\end{proof}

The following result generalises the result of Gower \cite{gower1982} to coEDMs and provides a new equivalent condition for a matrix to be coEDM.
\begin{corollary}\label{coEDM=wcnd}
    Let $D$ be a 
     hollow symmetric $n \times n$ matrix. Then $D$ is a conditional Euclidean distance matrix if and only if there exists $\omega \in e^\perp$ such that $y^TDy \leq 0$ for all $y \in {\rm span}\{1,\omega\}^\perp$.
\end{corollary}
\begin{proof}
    Suppose $D$ is a coEDM. Let $y \in {\rm span}\{1,\omega\}^\perp$. Then by Theorem \ref{proj_form}
    \[y^TDy = -2y^TFy \leq 0.\]
Conversely, suppose $y^TDy \leq 0$ for all $y \in {\rm span}\{1,\omega\}^\perp$. Let $s \in \rr^n$ such that $e^Ts=1$ and $s^T\omega=0$. Let $x \in \omega^\perp$
and take $y = (I-se^T)x$. Then $y = x-(e^Tx)s$ and hence $y^Te=y^T\omega=0$. Thus
\[ -\frac{1}{2}x^T(I-es^T)D(I-se^T)x = -\frac{1}{2}y^TDy \geq 0.\]
Therefore,  by Theorem \ref{proj_form},
$D$ is a coEDM. This completes the proof.
\end{proof}


    

The above result shows that $D$ is a coEDM if and only if $D$ is negative semidefinite on ${\rm span}\{1,\omega\}^\perp$. 
    For a given $\omega \in \rr^n$, we define  $H_{\omega}^n:=\{x \in \rr^n~:~e^Tx = \omega^Tx= 
    0\}$.
It is known that a non-zero EDM has exactly one positive eigenvalue. The following lemma shows that a coEDM can have at most two positive eigenvalues.


\begin{theorem}\label{two_pos_ev}
    Suppose $D$ is a  conditional Euclidean distance matrix. Then $D$ has at most two positive eigenvalues.
\end{theorem}
\begin{proof}
    Suppose $D$ is a coEDM. Then by Corollary \ref{coEDM=wcnd}, there exists $\omega \in e^\perp$ such that $y^TDy \leq 0$ for all $y \in H_{\omega}^n$. 
    
    If possible, suppose $D$ has three positive eigenvalues $\lambda$, $\mu$, and $\gamma$ with eigenvectors $x$, $y$ and $z$, respectively. Without loss of generality, we may assume that $x$, $y$ and $z$ are orthogonal. 

    Since $x^TDx = \lambda x^Tx>0$, $x \notin H_{\omega}^n$.  Similarly $y \notin H_{\omega}^n$ and $z\notin H_{\omega}^n$. Let
    \[c=(e^Tz)(\omega^Ty)-(e^Ty)(\omega^Tz),\]
    \[d=(e^Tz)(\omega^Tx)-(\omega^Tz)(e^Tx),\]
    and 
    \[f=(e^Ty)(\omega^Tx)-(\omega^Ty)(e^Tx).\]

In order to see that the constants $c$, $d$ and $f$ are non-zero, we need to check if the matrices 
\[\left[\begin{array}{cc}
    e^Tz & e^Ty  \\
     \omega^Tz& \omega^Ty
\end{array}\right], \left[\begin{array}{cc}
    e^Tz & e^Tx  \\
     \omega^Tz& \omega^Tx
\end{array}\right], \left[\begin{array}{cc}
    e^Ty & e^Tx  \\
     \omega^Ty& \omega^Tx
\end{array}\right]\]
are non-singular. It is enough to show that the columns of these matrices are linearly independent. If possible suppose there are non-zero constants $\alpha$ and $\beta$ such that 
\[\alpha \left[\begin{array}{cc}
    e^Tz  \\
     \omega^Tz
\end{array}\right] + \beta \left[\begin{array}{cc}
     e^Ty  \\
      \omega^Ty
\end{array}\right] = \0.\]
Then $\alpha z+ \beta y \in H_{\omega}^n$ and 
\begin{equation*}
    \begin{aligned}
        (\alpha z+ \beta y)^TD(\alpha z+ \beta y) = \alpha^2\gamma z^Tz+\beta^2 \mu y^Ty>0,
    \end{aligned}
\end{equation*}
which is a contradiction since $D$ is is negative semidefinite on $H_{\omega}^n$. Hence the constant $c$ is non-zero. It can be shown similarly that the constants $d$ and $f$ are non-zero.

Next, we note that  $cx-dy+ez \in H_{\omega}^n$ and 
\begin{equation*}
    \begin{aligned}
      (cx-dy+ez)^TD(cx-dy+ez)&=c^2x^TDx+d^2y^TDy+e^2z^TDz\\
      &= c^2\lambda x^Tx+d^2 \mu y^Ty+ e^2 \gamma z^Tz >0   ,
    \end{aligned}
\end{equation*}
which is again a contradiction. Hence $D$ has at most two positive eigenvalues. The proof is complete.
\end{proof}

One important aspect of EDMs is that they represent the distance between $n$ points in the Euclidean space. In our next result we obtain a similar representation of coEDMs and show that coEDMs represent perturbed distances between $n$ points in the Euclidean space. 
\begin{theorem}\label{coEDM_rep}
     A $n \times n$ 
      hollow symmetric matrix $D:=(d_{ij})$ is a conditional Euclidean distance matrix if and only if there exists vectors $p^1, p^2, \dotsc , p^n \in \rr^s$, for some $s$, 
     $\omega:=(\omega_i) \in \rr^n$, and  $\alpha \in \rr$, 
         such that
        \[d_{ij}=\|p^i-p^j\|^2+\alpha(\omega_i-\omega_j)^2 \]
        for $i,j=1,2,\dotsc,n$.
\end{theorem}
\begin{proof}
    Suppose $D$ is a coEDM.  From (\ref{coEDM_is_union}), we know that 
    \[\mathscr{C}^n=\mathscr{D}^n \cup \mathscr{K}(\Omega_n^{(1)}(e)). \]    
    If $D$ is an EDM, then the result follows by using (\ref{Def_EDM}).  Assume $D \in \mathscr{K}(\Omega_n^{(1)}(e))$.  Then there exists $F \in \Omega_n^{(1)}(e)$ such that $D=\mathscr{K}(F)$. Let $r={\rm rank}(F) \neq 1$ and $\lambda_1 \geq \dotsc \geq \lambda_{r-1}>0>\alpha$ be the ordered non-zero eigenvalues of $F$. We can write the spectral decomposition of $F$ as
\[F = \Gamma \Lambda \Gamma^T+ \alpha\omega\omega^T,\]
with $\Gamma = (\gamma_1,\dotsc, \gamma_{r-1})$, and $\Lambda = {\rm diag}(\lambda_1,\dotsc,\lambda_{r-1})_{(r-1) \times (r-1)}$, where $\Gamma^T\Gamma=I$, $\Gamma^Te=\0$, $\Gamma^T\omega =\0$ and $\omega^Te=0$.
Since $\Gamma\Lambda\Gamma^T \in \Omega_n(e)$, the matrix $\Delta=\mathscr{K}(\Gamma\Lambda\Gamma^T)$ is an EDM. Applying the map $\mathscr{K}$ to $F$, we get 
\[D=\mathscr{K}(F) = \Delta+\alpha \mathscr{K}(\omega\omega^T).\]
Let $\omega=(\omega_i)$. By (\ref{Def_EDM}), there exists vectors $p^1, p^2, \dotsc , p^n \in \rr^{r-1}$
         such that
        \[d_{ij}=\|p^i-p^j\|^2+\alpha(\omega_i-\omega_j)^2 \]
        for $i,j=1,2,\dotsc,n$.
If $r=1$, then the only non-zero eigenvalue of $F$ is negative. Therefore there exist $\alpha<0$ and $\omega \in \rr^n$ such that  $F=\alpha \omega \omega^T$. Choosing $p^1=p^2= \dotsc =p^n = \0 \in \rr^n$, we get
\[d_{ij}= \mathscr{K}(F) = \|p^i-p^j\|^2+\alpha(\omega_i-\omega_j)^2\]
 for $i,j=1,2,\dotsc,n$.
 
For the reverse implication, suppose there exist vectors $p^1, p^2, \dotsc , p^n \in \rr^s$, for some $s$, $\omega:=(\omega_i) \in \rr^n$ and $\alpha \in \rr$ such that
        \begin{equation}\label{coEDM_point_form}
            d_{ij}=\|p^i-p^j\|^2+\alpha (\omega_i-\omega_j)^2
        \end{equation}
        for $i,j=1,2,\dotsc,n$.
In (\ref{coEDM_point_form}), without loss of generality, we can assume  that $\alpha<0$, since if $\alpha>0$, the element $d_{ij}$ is written as $d_{ij}=\|y^i-y^j\|^2$ where $$
y^i=({p^i}^T, \sqrt{\alpha}\omega_i)^T \in \rr^{s+1},$$
which means $D$ is an EDM. Clearly, when $\alpha=0$, the matrix $D$ is also an EDM. Now suppose $\alpha<0$ in (\ref{coEDM_point_form}). 

We can also assume, without loss of generality, that $\omega^Te=0$. For, if $\omega^Te \neq0$, we can replace $\omega$ by $\omega - \frac{\omega^Te}{n}e$ in (\ref{coEDM_point_form}) without affecting any entry in $D$. Now, we have
        \begin{equation*}
            \begin{aligned}
                d_{ij}&=\|p^i\|^2+\|p^j\|^2-2\langle p^i, p^j \rangle+\alpha(\omega_i^2+\omega_j^2-2\omega_i\omega_j)\\
            &=\|p^i\|^2+\alpha\omega_i^2+\|p^j\|^2+\alpha\omega_j^2-2(\langle p^i, p^j \rangle+\alpha\omega_i\omega_j).
            \end{aligned}
        \end{equation*}
        Define $G:=(g_{ij})$, where $g_{ij}:= \langle p^i, p^j \rangle+\alpha\omega_i\omega_j$.  We claim that $G$ is positive semidefinite on $\omega^\perp$.
Note that 
\begin{equation*}
    \begin{aligned}
        G &= XX^T+ \alpha\omega \omega^T,
    \end{aligned}
\end{equation*}
where $X$ has $p^i$ as its $i^{\rm th}$ row, for each $i=1,\dotsc,n$. Suppose $x \in \omega^\perp$. Then
\begin{equation*}
    \begin{aligned}
        x^TGx &= x^TXX^Tx+ \alpha x^T\omega  \omega^Tx = \|x^TX\|^2 \geq 0.
    \end{aligned}
\end{equation*}
Thus $G$ is positive semidefinite on $\omega^\perp$.  It is easy to see that
        \[ D = E\Diag(G)+\Diag(G)E-2G.\]
Hence $-\frac{1}{2}JDJ=JGJ$. Define
\[F=-\frac{1}{2}JDJ.\]
Let $x \in \omega^\perp$. Then
\[x^TJ\omega = x^T(I-\frac{1}{n}E) \omega= x^T\omega = 0.\]
Thus $x^TJ \in \omega^\perp$ and therefore $x^T F x =x^TJ   GJx\geq 0$. Hence, $F$ is positive semidefinite on $\omega^\perp$. So, $D$ is a coEDM. This completes the proof.
\end{proof}
The above result shows that any $n$ points in a Euclidean space, together with an additional coordinate specified by $\omega$, can be used to generate a coEDM. The next result provides a necessary and sufficient condition for a coEDM to be non-negative and is a direct consequence of Theorem \ref{coEDM_rep}.


\begin{corollary}\label{coEDM_nonneg}
Suppose $D:=(d_{ij})$ is a conditional Euclidean distance matrix that is not an EDM and let $p^1, p^2, \dotsc , p^n \in \rr^s$, for some $s$,  $\omega:=(\omega_i) \in \rr^n$ and  $\alpha \in \rr$,  such that
        \[d_{ij}=\|p^i-p^j\|^2+\alpha(\omega_i-\omega_j)^2 \]
        for $i,j=1,2,\dotsc,n$.
Then $D$ is non-negative if and only if 
\[\alpha \geq -\underset{\{(i,j):\omega_i\neq \omega_j\}}{\rm min}\frac{\|p^i-p^j\|^2}{(\omega_i-\omega_j)^2}.\]
\end{corollary}

We will next illustrate the construction of coEDMs given in the proof of Theorem \ref{coEDM_rep} and the non-negativity condition given in Corollary \ref{coEDM_nonneg} through an example. The example uses the notations introduced in the proof of Theorem \ref{coEDM_rep}.
\begin{example} \rm
    Let 
    \[\gamma_1 = \frac{1}{\sqrt{2}}\left[\begin{array}{r}
         1  \\
         -1\\
         0\\
        0
    \end{array}\right], ~\gamma_2 = \frac{1}{\sqrt{2}}\left[\begin{array}{r}
         1  \\
         0\\
         -1\\
        0
    \end{array}\right],~\omega = \frac{1}{\sqrt{2}}\left[\begin{array}{r}
         1  \\
         0\\
         0\\
        -1
    \end{array}\right].\]
Define 
\[F_{\alpha}:=2\gamma_1\gamma_1^T+2\gamma_2\gamma_2^T+\alpha \omega\omega^T.\]
For each $\alpha \in \rr$, $F_{\alpha} \in \Omega_n^{-}(e)$ and hence $D_{\alpha}:=\mathscr{K}(F_{\alpha})$ is a coEDM. Since 
\[\mathscr{K}(vv^T) = ((v_i-v_j)^2)\]
for any vector $v$, we can easily compute $D_{\alpha}$ as follows:
\begin{equation*}
    \begin{aligned}
        D_{\alpha} &= 2\mathscr{K}(\gamma_1\gamma_1^T)+2\mathscr{K}(\gamma_2\gamma_2^T)+\alpha \mathscr{K}(\omega\omega^T)\\
        &=\left[\begin{array}{cccc}
            0 & 4 &1 & 1 \\
             4 & 0 &1 & 1 \\
             1 & 1 &0 & 0 \\
             1 & 1 &0 & 0 
        \end{array}\right]+\left[\begin{array}{cccc}
            0 & 1 &4 & 1 \\
             1 & 0 &1 & 0 \\
             4 & 1 &0 & 1 \\
             1 & 0 &1 & 0 
        \end{array}\right]+\frac{\alpha}{2}\left[\begin{array}{cccc}
            0 & 1 &1 & 4 \\
             1 & 0 &0 & 1 \\
             1 & 0 &0 & 1 \\
             4 & 1 &1 & 0 
        \end{array}\right]\\
        &= \left[\begin{array}{cccc}
            0 & 5 &5 & 2 \\
             5 & 0 &2 &1 \\
             5 & 2 &0 & 1 \\
             2 & 1 &1 & 0 
        \end{array}\right]+\frac{\alpha}{2}\left[\begin{array}{cccc}
            0 & 1 &1 & 4 \\
             1 & 0 &0 & 1 \\
             1 & 0 &0 & 1 \\
             4 & 1 &1 & 0 
        \end{array}\right]=(\|p^i-p^j\|^2+\alpha(\omega_i-\omega_j)^2).
    \end{aligned}
\end{equation*}
Since 
    \[-\underset{\{(i,j):\omega_i\neq \omega_j\}}{\rm min}\frac{\|p^i-p^j\|^2}{(\omega_i-\omega_j)^2} = -\frac{\|p^1-p^4\|^2}{(\omega_1-\omega_4)^2}=-1,\]
the coEDM, $D_{\alpha}$ is non-negative if and only if $\alpha \geq -1$. Taking $\alpha =-\frac{1}{2}$, we obtain the following non-negative coEDM.
\begin{equation*}
    \begin{aligned}
        D_{\alpha} &= \left[\begin{array}{cccc}
            0 & 5 &5 & 2 \\
             5 & 0 &2 &1 \\
             5 & 2 &0 & 1 \\
             2 & 1 &1 & 0 
        \end{array}\right]-\frac{1}{4}\left[\begin{array}{cccc}
            0 & 1 &1 & 4 \\
             1 & 0 &0 & 1 \\
             1 & 0 &0 & 1 \\
             4 & 1 &1 & 0 
        \end{array}\right]= \frac{1}{4}\left[\begin{array}{cccc}
            0 & 19 &19 & 4 \\
             19 & 0 &8 & 3 \\
             19 & 8 &0 & 3 \\
             4 & 3 &3 & 0 
        \end{array}\right].
    \end{aligned}
\end{equation*}
We further note that for $\alpha=-\frac{1}{2}$, $D_{\alpha}$ is not an EDM since $F_{\alpha}$ has one negative eigenvalue.
\end{example}

\subsection{Projected centered symmetric matrix}
We recall that $J=I-\frac{1}{n}E$ is the orthogonal projection matrix onto $e^\perp$. A different choice of basis of  $e^\perp$ will give rise to different factorizations of $J$. We present here two such factorizations. Define 
$U_{n \times (n-1)}:=\left[\begin{array}{cc}
     -e^T  \\
     I
     \end{array}\right]$
and $V$ to be a $n \times (n-1)$ matrix whose columns form an orthonormal basis of $e^\perp$. It is easy to see that $J = UU^\dagger = VV^T$. The matrix $V$ is not unique. We choose the following $V$
\[V = \left[\begin{array}{cc}
     ye^T \\
     I+xee^T
\end{array}\right],~{\rm where}~ y = -\frac{1}{\sqrt{n}}~{\rm and}~x=-\frac{1}{n+\sqrt{n}}.\]
We note that 
\[V=U(I+xee^T), ~{\rm and }~U=V(I+\frac{x}{y}ee^T).\]
Consider the following two transformations $\mathscr{T}_V: \mathscr{S}^n_h \rightarrow \mathscr{S}^{n-1}$ and $\mathscr{K}_V: \mathscr{S}^{n-1} \rightarrow \mathscr{S}^n_h$ defined by
\[ \mathscr{T}_V(A) = V^T\mathscr{T}(A)V = -\frac{1}{2}V^TAV,~\mbox{and}\]
\[\mathscr{K}_V(X) = \mathscr{K}(VXV^T)= E\Diag(VXV^T)+\Diag(VXV^T)E-2VXV^T.\]
We call $\mathscr{T}_V(A)$ the \textit{projected centered symmetric matrix} of $A$. It is easy to see that $\rank(\mathscr{T}(A)) = \rank(\mathscr{T}_V(A))$ and $\mathscr{T}(A)=V\mathscr{T}_V(A)V^T$ for all $A \in \mathscr{S}_h^n$. Furthermore, the transformation $\mathscr{T}_V$ restricted to $\mathscr{S}_h^n$ is the inverse of the transformation $\mathscr{K}_V$ (see \cite{alfakih}, Chapter 4).  

\begin{theorem}\label{thm_proj_cen}
    Let $D$ be an $n \times n$ hollow symmetric matrix. Then $D$ is a coEDM if and only if $\mathscr{T}_V(D)$ is positive semidefinite on $\{\tilde{\omega}\}^\perp$ for some $\tilde{\omega} \in \rr^{n-1}$.
\end{theorem}
\begin{proof}
    Using Theorem \ref{def_coEDM}, it is enough to show that $\mathscr{T}_V(D)$ is positive semidefinite on $\tilde{\omega}^\perp$ for some $\tilde{\omega} \in \rr^{n-1}$ if and only if $\mathscr{T}(D)$ is positive semidefinite on $\omega^\perp$ for some $\omega \in e^\perp$.

Let $\omega \in e^\perp$ such that $\mathscr{T}(D)$ is positive semidefinite on $\omega^\perp$. Define $\tilde{\omega}:= V^T\omega \in \rr^{n-1}$ and let $x \in \tilde{\omega}^\perp$. Since
\[x^T\mathscr{T}_V(D)x = x^TV^T\mathscr{T}(D)Vx,\]
and $x^TV^T\omega = x^T\tilde{\omega} = 0$, it follows that $x^T\mathscr{T}_V(D)x \geq 0$. Conversely, suppose $\tilde{\omega} \in \rr^{n-1}$ such that $\mathscr{T}_V(D)$ is positive semidefinite on $\{\tilde{\omega}\}^\perp$. Define $\omega:= V\tilde{\omega} \in \rr^{n}$.
Clearly $\omega \in e^\perp$. Let $x \in \omega^\perp$, and consider
\[x^T\mathscr{T}(D)x = x^TV\mathscr{T}_V(D)V^Tx.\]
Since $x^TV \tilde{\omega} = x^T\omega=0$, it is evident that $x^T\mathscr{T}(D)x \geq 0$. The proof is complete.
\end{proof}

\begin{theorem}
    Let $D = \left[\begin{array}{cc}
         0& d^T  \\
         d& \tilde{D} 
    \end{array}\right]$ be an $n \times n$ real symmetric matrix whose diagonal entries are all zero, where $d \in \rr^{n-1}$. Then $D$ is an coEDM if and only if 
    \[de^T+ed^T-\tilde{D}\]
    is positive semidefinite on $\omega^\perp$ for some $\omega \in \rr^{n-1}$.
\end{theorem}
\begin{proof}
From Theorem  \ref{thm_proj_cen}, $D$ is a coEDM if and only if $-V^TDV$ is positive semidefinite on $\omega^\perp$ for some $\omega \in \rr^{n-1}$. But $V=UA$, where $A$ is the non-singular matrix $I+xee^T$. Thus $D$ is a coEDM if and only if   $-U^TDU$ is positive semidefinite on $\omega^\perp$ for some $\omega \in \rr^{n-1}$. The result now follows because $U^TDU = -de^T-ed^T+\tilde{D}$. 
\end{proof}

\section{The pre-Gale space and the pre-Gale matrix}
Let $D$ be an $n \times n$ EDM. Then for some positive integer $s$, there exist points $p^1, \dotsc , p^n \in \rr^s$ such that
\begin{equation*}
    d_{ij}= \|p^i-p^j\|^2
\end{equation*}
for each $i, j \in \{1,\dotsc,n\}$. The dimension of the affine span of these $n$ points is called the \textit{embedding dimension} of $D$. Let $X:=[p^1,\dotsc,p^n]^T$ and $F=XX^T$. It is easy to see that $F=-\frac{1}{2}JDJ$. The \textit{Gale space} of an
$n \times n$ EDM $D$ of embedding dimension $r$, denoted by ${\rm gal}(D)$, is defined as
\[{\rm gal}(D) = {\rm null}(X^T) \cap {\rm null}(e^T). \]
It is well known that for EDMs with embedding dimension $r$, the dimension of the Gale space is $n-1-r$. Since ${\rm null}(F) = {\rm null}(X^T)$, Gale space can also be defined as 
\[{\rm gal}(D) = {\rm null}(F) \cap {\rm null}(e^T). \]
Any $n \times (n-1-¯r)$ matrix $Z$ whose columns form a basis of ${\rm gal}(D)$ is called a \textit{Gale matrix} of $D$. For EDMs, the Gale space and the Gale matrix have been studied in detail.   The geometric significance of the Gale space and the Gale matrix is also studied in detail in the literature (see Alfakih \cite{alfakih}).
Motivated by the significance of the Gale space in the theory of EDMs, we define the following for a hollow symmetric matrix.

\begin{definition}\rm
The \textit{pre-Gale space} of an $n \times n$ hollow symmetric matrix $D$, denoted by $\mbox{pGal}(D)$, is defined as
\[\mbox{pGal}(D) = {\rm null}(F) \cap {\rm null}(e^T), \]
where $F=-\frac{1}{2}JDJ$.    
\end{definition}

 Note that if $D$ is an EDM then its pre-Gale space and Gale space coincide. Let $D$ be an $n \times n$ hollow symmetric matrix and $\tilde{r}$ be the dimension of $\mbox{pGal}(D)$. Clearly $\tilde{r}= n-1-\rank(F)$. 
 
 \begin{definition}\rm 
     Any $n \times \tilde{r}$ matrix $Z$ whose columns form a basis of $\mbox{pGal}(D)$ is said to be a \textit{pre-Gale matrix} of $D$. 
 \end{definition}

We now present some results on the pre-Gale space and the pre-Gale matrix of a hollow symmetric matrix. Our first result establishes a connection between pre-Gale matrices and a  hollow symmetric matrix. 

\begin{lemma}
    Let $D$ be a hollow symmetric matrix and let $Z$ be a pre-Gale matrix of $D$. Further, let $F = -\frac{1}{2}JDJ$, then
    \[DZ = e \xi^T,\]
    where $\xi = Z^T\diag(F).$
\end{lemma}
\begin{proof}
    Since the maps $\mathscr{T}$ and $\mathscr{K}$ are mutually inverse, we have
    \[D = E\Diag(F)+\Diag(F)E-2F.\]
    The result now follows from the fact that the columns of $Z$ form a basis for $\mbox{pGal}(D) = {\rm null}(F) \cap {\rm null}(e^T)$. 
\end{proof}

The following result relates the affine dependence of columns of $F$ to the linear dependence of columns of a {pre-Gale matrix} $Z$.

\begin{lemma}\label{lem_affine_dep}
    Let $\rank(F)=r$, $f^1, \dotsc,f^n$ be columns of $F$ and $z^1, \dotsc,z^n \in \rr^{n-r-1}$ be the rows of a {pre-Gale matrix} $Z$. Let $\mathcal{S}$ be a subset of $\{1,\dotsc,n\}$ of cardinality $r+1$. Then $\{f^i:i \in \mathcal{S}\}$ are affinely dependent if and only if $\{z^i:i \in \mathcal{S}^c\}$ are linearly dependent.
\end{lemma}
\begin{proof}
    Without loss of generality, we assume that $\mathcal{S}=\{1,\dotsc,r+1\}$ i.e. assume that $f^1,\dotsc,f^{r+1}$ are affinely dependent. Then there exists a non-zero $\lambda=(\lambda_i) \in \rr^{r+1}$ such that
    \[\sum_{i=1}^{r+1} \lambda_if^i=\0~{\rm and}~\sum_{i=1}^{r+1}\lambda_i=0.\]
    Let $x=\left[\begin{array}{cc}
        \lambda^T & \0^T  
    \end{array}\right]^T \in \rr^n$. Then
    \[\sum_{i=1}^{n} x_if^i=\0~{\rm and}~\sum_{i=1}^{n}x_i=0.\]
Thus $x \in \mbox{pGal}(D)$. Partition the pre-Gale matrix $Z$ as $Z=\left[\begin{array}{cc}
     Z_1  \\
     Z_2
\end{array}\right]$, where $Z_1$ is of order $(r+1) \times (n-r-1)$. Then
\[\left[\begin{array}{cc}
     Z_1  \\
     Z_2
\end{array}\right]y=x=\left[\begin{array}{cc}
        \lambda \\ \0  
    \end{array}\right]\]
for some non-zero $y \in \rr^{n-r-1}$. Consequently, the matrix $Z_2$ is singular. Therefore the rows $z^{r+2},\dotsc,z^n$ are linearly dependent. The proof of the converse follows by reversing the above steps. This completes the proof.
\end{proof}

\begin{corollary}
     Let $\rank(F)=r$ and $z^1, \dotsc,z^n \in \rr^{n-r-1}$ be the rows of a {pre-Gale matrix} $Z$. If any $r+1$ columns of $F$ are affinely independent, then there exists a pre-Gale matrix of the form $Z=\left[\begin{array}{cc}
     \tilde{Z}  \\
    I_{n-r-1}
\end{array}\right]$.
\end{corollary}
\begin{proof}
    Without loss of generality, assume that $f^1,\dotsc,f^{r+1}$ are affinely independent. Partition the pre-Gale matrix $Z$ as $Z=\left[\begin{array}{cc}
     Z_1  \\
     Z_2
\end{array}\right]$, where $Z_1$ is of order $(r+1) \times (n-r-1)$. Then by Lemma \ref{lem_affine_dep}, $Z_2$ is non-singular. The result follows by multiplying $Z$ on the right with $Z_2^{-1}$. 
\end{proof}

\section{Non-negative coEDMs}
In this last section we obtain some results for non-negative coEDMS using their pre-Gale space. These results are generalisations of results known for EDMs (see Alfakih \cite{alfakih}, Chapter 3). We establish a connection between the null space of a non-negative coEDM and its pre-Gale space in our next result. 
\begin{theorem}\label{nullD_subset_galD}
    Let $D$ be a non-zero, non-negative coEDM of order $n$ and let $\omega \in \rr^n$ such that $F = -\frac{1}{2}JDJ$ is positive semidefinite on $\omega^\perp$. Then
    \[{\rm null}(D) \cap {\rm null}(\omega^T) \subseteq \mbox{pGal}(D).\]
\end{theorem}
\begin{proof}
    Let $x \in {\rm null}(D) \cap {\rm null}(\omega^T)$ and let $F = -\frac{1}{2}JDJ$. Then
    \[2x^TFx = -\frac{(e^Tx)^2 (e^TDe)}{n^2}.\]
    Since $F$ is positive semidefinite on $\omega^\perp$, $x^TFx \geq 0$. Also, $e^TDe >0$. Thus, it follows that $ e^Tx = 0$ and $x^TFx = 0$.  This proves that $ x \in \mbox{pGal}(D)$ and hence the result follows.
\end{proof}

The next result provides a sufficient condition for $e$ to be in the column space of a coEDM. 

\begin{theorem}\label{1_in_colD}
    Let $D$ be a non-zero, non-negative coEDM of order $n$ and let $\omega \in \rr^n$ such that $F = -\frac{1}{2}JDJ$ is positive semidefinite on $\omega^\perp$. If ${\rm null}(D) \subseteq \omega^\perp$, then $e \in {\rm col}(D)$.
\end{theorem}
\begin{proof}
    From Theorem \ref{nullD_subset_galD}, we have
    \[\mbox{pGal}(D)^\perp \subseteq {\rm col}(D).\]
    The result follows from the fact that $e \in \mbox{pGal}(D)^\perp $.
\end{proof}

\begin{theorem}\label{1_in_colD2}
    Let $D$ be a non-zero, non-negative coEDM of order $n$ and let $\omega \in \rr^n$ such that $F = -\frac{1}{2}JDJ$ is positive semidefinite on $\omega^\perp$. If ${\rm null}(D) \subseteq \omega^\perp$, then 
    \[\rank(D) = \rank(F)+1 \mbox{ \rm or }\rank(D) = \rank(F)+2.\]
\end{theorem}
\begin{proof}
  Since 
  \[  D=E\Diag(F)+\Diag(F)E-2F,\]
  $\rank(D) \leq \rank(F)+2$. Further, $\mbox{pGal}(D)^\perp \subseteq {\rm col}(D)$ implies $\rank(D) \geq \rank(F)+1$, since ${\rm dim}(\mbox{pGal}(D))= n-1-\rank(F)$. The proof is complete.
\end{proof}

\section{Future scope: connection to Minkowski spaces, machine learning, and general relativity}

Euclidean spaces play an important role for the representation of data in machine learning.  Representing data in a Euclidean space is appropriate to compare lengths and angles. However, Euclidean space often fail to capture some special structures like hierarchical structures in the data. Since hierarchies can be formulated as trees, hyperbolic spaces are used in literature to represent hierarchically structured data (see \cite{gromov1987hyperbolic}).  In \cite{pmlr-v97-law19a}, Marc et al. proposed the representation of hierarchical data using the squared Lorentzian distance and showed it is a better choice as compared to the other commonly used metrics. 

\begin{definition}
The $(n + 1)$-dimensional
Minkowski space $\mathbb{M}^{1,n}$ is the real  vector space $\rr^{n+1}$
equipped with the symmetric bilinear form:    
\[ \langle x, y \rangle_{\mathcal{L}} = -x_0y_0+x_1y_1+\dotsc+x_ny_n,\]
where $x=(x_0,x_1,\dotsc,x_n)^T$ and $y=(y_0,y_1,\dotsc,y_n)^T \in \rr^{n+1}$. The symmetric bilinear form $\langle . , . \rangle_{\mathcal{L}}$ is called the \textit{Lorentzian inner product}. Further, the Lorentzian norm of a vector $x$ is defined to be 
\[\|x\|_\mathcal{L} = \sqrt{\langle x, x \rangle_\mathcal{L}}.\]
For $x,y \in \mathbb{M}^{1,n}$, $\|x-y\|_\mathcal{L}^2$ is called the squared Lorentzian distance between $x$ and $y$.
\end{definition}
The literature on Lorentzian inner product can be found in \cite{fan-etal-2024-enhancing} and \cite{Ratcliffe}.
Suppose $D=(d_{ij})$ is a coEDM that is not an EDM. Then by Theorem \ref{coEDM_rep}, there exists $p^1,\dotsc,p^n \in \rr^s$ for some $s$, $\omega:=(\omega_i) \in \rr^n$ and $\alpha<0$ such that 
\[d_{ij}= \|p^i-p^j\|^2+\alpha (\omega_i-\omega_j)^2\]
for $i,j =1,\dotsc,n$. This means
\[d_{ij}= \sum_{k=1}^{s}(p^i_k-p^j_k)^2-(\sqrt{-\alpha}\omega_i-\sqrt{-\alpha}\omega_j)^2\]
for $i,j =1,\dotsc,n$.
Define for each $i=1,\dotsc,n$, 
\[y^i = (\sqrt{-\alpha}\omega_i,p^i_1,\dotsc,p^i_s) \in \rr^{s+1}.\]
Then it is easy to see that 
\begin{equation*}
    \begin{aligned}
        d_{ij}&= \langle y^i-y^j, y^i-y^j \rangle_\mathcal{L}\\
        &= \| y^i-y^j\|_\mathcal{L}^2.
    \end{aligned}
\end{equation*}
for $i,j =1,\dotsc,n$. This shows that the entries of an $n \times n$ coEDM (which is not an EDM) represent the pairwise squared Lorentzian distance between $n$ points in a Minkowski space. 

\begin{definition}
    A subset $S \subseteq \mathbb{M}^{1,n}$ is called achronal if for every pair of distinct points in $x,y \in S$ 
    \[\|x-y\|_\mathcal{L}^2 \geq 0.\]
\end{definition}
The above definition of an achronal set can be easily deduced by using the definition of $1$-Lipschitz function and the fact that  every
achronal set is precisely the graph of a $1$-Lipschitz function (see \cite{BarbotMerigot2012}). Extensive literature on achronal sets, their physical interpretation and uses can be found in the book titled General relativity by Robert Wald \cite{Wald1984}.   

The next result is a direct consequence of the above discussion.

\begin{theorem}
    Suppose $D:=(d_{ij})$ is a $n \times n$ coEDM that is not an EDM and let $y^1,\dotsc,y^n \in \mathbb{M}^{1,s}$ such that 
    \[d_{ij}=\| y^i-y^j\|_\mathcal{L}^2.\]
    Then $D$ is non-negative if and only if the set $\{y^1,\dotsc,y^n \}$ is an achronal subset of $\mathbb{M}^{1,s}$.
\end{theorem}

There are several future directions for the continuation of theory of coEDMs and its applications. First, 
the interlink between coEDMs and distances in Minkowski spaces can be used to make advances in the theory of coEDMs from a geometric point of view. Second, the existing use of squared Lorentzian distance in the machine learning algorithm suggests coEDMs can be a useful mathematical tool for the representation of hierarchical data. Lastly, achronal sets allow physicists to track how information can and cannot flow through a curved universe. This connection to the theory of general relativity suggest the use of mathematical theory of coEDMs in physics.
\bibliography{mybibfile}

\noindent {\em Authors' address}:\\
{\em Shivani Goel}, School of Analytics, Tata Institute of Social Sciences, Mumbai, India.
 E-mail: \texttt{shivani.goel@tiss.ac.in} \\
\noindent 
{\em Hiroshi Kurata}, Faculty of Economics, Kanagawa University, Kanagawa, Japan.
 E-mail: \texttt{ft102272aj@kanagawa-u.ac.jp}

\end{document}